\documentclass[12pt, reqno]{amsart}
\usepackage[utf8]{inputenc}
\usepackage{amsmath}
\usepackage{amsfonts}
\usepackage{amssymb}
\usepackage{amsthm}
\usepackage{bm}
\usepackage{bbm}
\usepackage{xcolor}
\usepackage{url}
\usepackage[hidelinks]{hyperref}

\usepackage{graphicx,tikz}

\numberwithin{equation}{section}

\theoremstyle{plain}
\newtheorem{theorem}{Theorem}[section]
\newtheorem{lemma}[theorem]{Lemma}

\newtheorem{proposition}[theorem]{Proposition}

\newtheorem{conjecture}[theorem]{Conjecture}
\newtheorem{remark}[theorem]{Remark}

\title[Powers with small exponent from products in progression]{Perfect powers with small exponents from products of terms in progression}
\author{Kyle Pratt}
\address{Brigham Young University, Department of Mathematics, Provo, UT 84602, USA}
\email{kyle.pratt@mathematics.byu.edu}


\allowdisplaybreaks

\begin{document}
\date{}

\begin{abstract}
Let $k$ be a large integer, and let $\ell$ be an integer with $2 \leq \ell \leq (\log \log k)^{1/5}$. We determine all solutions to the Diophantine equation
\begin{align*}
n(n+d)\cdots (n+(k-1)d)=y^\ell
\end{align*}
in nonzero integers $n,d,y$ with $d\geq 1$ and $\text{gcd}(n,d)=1$. We also prove results for generalizations where the product $n(n+d)\cdots (n+(k-1)d)$ is ``almost'' a perfect power, rather than exactly a perfect power, and where some of the terms $n+id$ are dropped from the product.
\end{abstract}

\maketitle

\section{Introduction}

\subsection{History and main results}

Erd\H{o}s and Selfridge \cite{ES1975} proved that the product of two or more consecutive positive integers is never a perfect power. That is, they showed
\begin{align}\label{eq:ES eqn}
\prod_{i=0}^{k-1}(n+i) = y^\ell
\end{align}
has no solutions in positive integers $n,y,k,\ell$, where $k$ and $\ell$ are both at least two. The earliest results on this problem go back to Liouville \cite{Liouville1857} in 1857, who proved $2 \cdot 3 \cdots n$ is never a perfect power. 

Since the work of Erd\H{o}s and Selfridge, a great deal of work has been done on generalizations of \eqref{eq:ES eqn}. One way to generalize the situation is to take products of consecutive terms in a primitive arithmetic progression. Hence, one asks whether there are solutions to
\begin{align}\label{eq:eqn of study, no b}
\prod_{i=0}^{k-1}(n+id) = y^\ell
\end{align}
in positive integers $n,d,y,k,\ell$ with $\text{gcd}(n,d)=1$ and $k,\ell \geq 2$. There are infinitely many solutions to \eqref{eq:eqn of study, no b} if $k=2$ and if $(k,\ell)=(3,2)$, but otherwise one expects no solutions. More modestly, one may conjecture that \eqref{eq:eqn of study, no b} has no solutions once $k$ is sufficiently large. These conjectures are often attributed to Erd\H{o}s (see, e.g., \cite{BloomErdos672} and \cite[pp. 53-54]{MPE2013}, although Fan Chung notes this chapter was actually written by Ron Graham \cite{FCG}), and are known if one attaches additional conditions, like if $\ell$ and the number of prime factors of $d$ are fixed (work of Shorey and Tijdeman \cite{ST1990}), or if $n$ is fixed and $\ell \geq 7$ (work of Shorey \cite{Sho1996}). There are also results on solutions to \eqref{eq:eqn of study, no b} when $k$ is small (see \cite{BBGH2006,Bennett2018,GHP2009,Oblath1933}). The ABC conjecture implies that any solution to \eqref{eq:eqn of study, no b} must have $k$ bounded by an absolute constant (an unpublished argument due to Granville, but a proof later appeared in \cite[Section 6]{LS2015}). The recent book \cite{Sar2025} is entirely devoted to studying \eqref{eq:eqn of study, no b}, and we invite the reader to look there for further results and details.

Bennett and Siksek \cite{BS2020} made a spectacular breakthrough in the study of \eqref{eq:eqn of study, no b}. They used sophisticated techniques relating to the modularity of elliptic curves to show that, if $k$ is sufficiently large, then \eqref{eq:eqn of study, no b} has no solutions with $\ell >  \exp(10^k)$ prime \cite[Theorem 2]{BS2020}. If $k$ is large and fixed, this reduces the study of \eqref{eq:eqn of study, no b} with prime $\ell$ to a finite number of cases. Work of Darmon and Granville \cite{DG1995} then implies that \eqref{eq:eqn of study, no b} has at most finitely many solutions for fixed $k$.

Bennett and Siksek actually show that \eqref{eq:eqn of study, no b} has no solutions if $k$ is sufficiently large, $n$ and $y$ are nonzero integers, $\text{gcd}(n,d)=1$, and $\ell > \exp(10^k)$ is prime. The relaxation to nonzero integers $n$ and $y$ is significant, because \eqref{eq:eqn of study, no b} can then have solutions in some cases. For instance, when $k=4$ we take $n=-3$ and $d=2$ to get
\begin{align*}
(-3)(-1)(1)(3) = 3^2,
\end{align*}
and when $k=9$ we take $n=-7$ and $d=2$ to get
\begin{align*}
(-7)(-5)(-3)(-1)(1)(3)(5)(7)(9) = 315^2.
\end{align*}
When $(k,\ell) = (4,3)$ we can take $n \in \{-9,-6\}$ and $d=5$ to get
\begin{align*}
(-6)(-1)(4)(9) = (-9)(-4)(1)(6)=6^3.
\end{align*}
In general, we have the following conjecture, which is a stronger version of the conjectures attributed to Erd\H{o}s.

\begin{conjecture}\label{conj:main conj}
Let $k \geq 4$ and $\ell \geq 2$ be integers. If there are nonzero integers $n,d,y$ with $d \geq 1$ and $\textup{gcd}(n,d)=1$ such that \eqref{eq:eqn of study, no b} holds, then
\begin{itemize}
\item $(k,\ell,n,d) = (4,3,-9,5)$ or $(4,3,-6,5)$, or
\item $4 \mid k$ and $(\ell,n,d)=(2,-k+1,2)$, or
\item $k$ is an odd square and $(\ell,n,d)=(2,-k+2,2)$.
\end{itemize}
\end{conjecture}

The restriction to $k \geq 4$ in Conjecture \ref{conj:main conj} is necessary for finiteness, since there are infinitely many solutions to \eqref{eq:eqn of study, no b} when $(k,\ell)=(2,2),(3,2)$.

Some of the results mentioned above prove Conjecture \ref{conj:main conj} for small values of $k$, but in general it seems quite difficult. For example, Conjecture \ref{conj:main conj} implies a conjecture of Sander (see \cite[p. 2250]{BS2016} and \cite[Conjecture A]{Pratt2024}) regarding the rational points on the Erd\H{o}s--Selfridge superelliptic curves $y^\ell = x(x+1)\cdots (x+k-1)$. The problem of finding rational points on these curves is another way to generalize \eqref{eq:ES eqn}, and is, in some sense, of intermediate difficulty between \eqref{eq:ES eqn} and \eqref{eq:eqn of study, no b}.

The work of Bennett and Siksek resolves Conjecture \ref{conj:main conj} when $k$ is large and $\ell > \exp(10^k)$ is prime; thus, in the regime they study, $\ell$ is extremely large compared to $k$. In this work we are interested in the opposite regime, where $k$ is large and $\ell$ is very small compared to $k$. This regime is both interesting and complicated in light of Conjecture \ref{conj:main conj}, since \eqref{eq:eqn of study, no b} is conjectured to have solutions in some cases, and no solutions in other cases.

We show that, when $k$ is large and $\ell$ is small, the only solutions to \eqref{eq:eqn of study, no b} with nonzero $n$ and $y$ are those predicted by Conjecture \ref{conj:main conj}.

\begin{theorem}\label{thm:explicit solutions}
Let $k_0$ be a sufficiently large positive constant, and let $k \geq k_0$ be an integer. Let $\ell$ be an integer with $2 \leq \ell \leq (\log \log k)^{1/5}$. If there are nonzero integers $n,d,y$ with $d \geq 1$ and $\textup{gcd}(n,d)=1$ such that \eqref{eq:eqn of study, no b} holds, then 
\begin{itemize}
\item $4 \mid k$ and $(\ell,n,d)=(2,-k+1,2)$, or
\item $k$ is an odd square and $(\ell,n,d)=(2,-k+2,2)$.
\end{itemize}
\end{theorem}

By Theorem \ref{thm:explicit solutions} and the work of Bennett and Siksek \cite[Theorem 2]{BS2020}, if there is a counterexample to Conjecture \ref{conj:main conj} for some large $k$ and prime $\ell$, then we must have
\begin{align*}
(\log \log k)^{1/5} < \ell \leq \exp(10^k).
\end{align*}
The symmetry between the doubly-logarithmic lower bound and the doubly-exponential upper bound is amusing, but coincidental.

In fact, the upper bound on $\ell$ in Theorem \ref{thm:explicit solutions} is quite restrictive, and it would be desirable to have a result in which $\ell$ is allowed to be larger with respect to $k$. The condition on $\ell$ is imposed by Lemma \ref{lem:remonds theorem} below (see Remarks \ref{rmk:upper bound on ell not optimal given remond} and \ref{rmk:better bds on num of rat pts} for additional discussion). It seems, however, that our method does not allow one to take $\ell \gtrapprox \log k$.

The constant $k_0$ in Theorem \ref{thm:explicit solutions} is effectively computable. Hence, for a fixed value of $\ell$, resolving Conjecture \ref{conj:main conj} for that choice of $\ell$ could, in principle, be reduced to solving \eqref{eq:eqn of study, no b} for a finite number of values of $k$. However, the value of $k_0$ produced by our method is likely extremely large. It would be interesting to see if our method could be refined or modified in order to yield a respectable result for $\ell=2$, say.

We note that Theorem \ref{thm:explicit solutions} proves Sander's conjecture on rational points on Erd\H{o}s--Selfridge superelliptic curves (mentioned above) when $k$ is large and $2\leq \ell \leq (\log \log k)^{1/5}$. We leave the details to the interested reader.

Many of the works studying \eqref{eq:eqn of study, no b} actually apply to a more general equation. For a nonzero integer $b$, let $P^+(b)$ denote the largest prime dividing $b$ if $|b| >1$, otherwise define $P^+(b)=1$. The more general question is to ask for solutions to
\begin{align}\label{eq:eqn of study, all terms, by^ell}
\prod_{i=0}^{k-1}(n+id) = by^\ell, \ \ \ \ \ P^+(b) \leq k.
\end{align}
Thus, any prime $>k$ that divides the left-hand side must appear with a multiplicity divisible by $\ell$, but there are no restrictions on primes $\leq k$.

Note that one can indeed have solutions to \eqref{eq:eqn of study, all terms, by^ell} with $d=1$ or $d=2$. The existence of solutions for $d=2$ is noted above. When $d=1$, one can take $n=1$, say, and we have
\begin{align*}
\prod_{i=0}^{k-1}(n+i) = k! = k! \cdot 1^\ell,
\end{align*}
so we have a solution to \eqref{eq:eqn of study, all terms, by^ell} with $b=k!$ and any integer $\ell \geq 2$.

In fact, \eqref{eq:eqn of study, all terms, by^ell} often has many solutions. For any $k\geq 5$, let $p = p(k)$ be the least prime $>k$. Each of the following choices of $(n,d)$ then gives a solution to \eqref{eq:eqn of study, all terms, by^ell} for the given value of $k$:
\begin{equation}\label{eq:solutions with d=1,2 and b present}
\begin{alignedat}{4}
&(n,1), \quad &1      \le &n \le p-k,\\
&(n,1),&-(p-1) \le &n \le -k,\\
&(n,2),&-(p-2) \le &n \le p-2-2(k-1),\qquad n\text{ odd}.
\end{alignedat}
\end{equation}
Thus, \eqref{eq:eqn of study, all terms, by^ell} always has at least $3(p(k)-k)$ admissible pairs $(n,d)$. Note that every solution in \eqref{eq:solutions with d=1,2 and b present} leads to $P^+ \left( \prod_{i=0}^{k-1}(n+id)\right) \leq k$, so no primes $>k$ appear in the factorization of the left-hand side of \eqref{eq:eqn of study, all terms, by^ell}.

When $k$ is small, \eqref{eq:solutions with d=1,2 and b present} does not exhaust the pairs $(n,d)$. We have the additional pairs $(n,d)=(-16,7),(-12,7),(-8,3),(-4,3)$ when $k=5$, and the pairs $(n,d)=(-10,3),(-8,3)$ when $k=7$. These additional pairs still have the property that no term $n+id$ is divisible by a prime $>k$. In general, we expect the following to be true.

\begin{conjecture}\label{conj:b factor present}
Let $k \geq 5, \ell \geq 2$ be integers. If there are nonzero integers $n,d,b,y$ with $d \geq 1$, $\textup{gcd}(n,d)=1$, and $P^+(b) \leq k$ such that \eqref{eq:eqn of study, all terms, by^ell} holds, then $P^+ \left( \prod_{i=0}^{k-1}(n+id)\right) \leq k$. If $k \geq 8$, then $d \in \{1,2\}$, and $(n,d)$ is one of the pairs in \eqref{eq:solutions with d=1,2 and b present}.
\end{conjecture}

The restriction to $k \geq 5$ is again for finiteness, since $n(n+d)(n+2d)(n+3d)=6y^2$ has infinitely many solutions (a related elliptic curve has positive rank).

According to Conjecture \ref{conj:b factor present}, the equation \eqref{eq:eqn of study, all terms, by^ell} has $3(p(k)-k)$ admissible pairs $(n,d)$ when $k \geq 8$. It is difficult to give precise upper bounds on the number of solutions, in terms of $k$, since such a bound depends on results on prime gaps. If $n$ is positive, say, then $n$ must be quite small in terms of $k$, otherwise one can use results on primes in short intervals to find a term $n+id$ that is a prime $>k$, which would contradict the conclusion of the conjecture. In particular, work of Baker, Harman, and Pintz \cite{BHP2001} would show $n \lessapprox k^{0.525}$. Assuming some version of a Cramer conjecture about the size of gaps between consecutive primes (see \cite{BFT2023} and \cite{GL2023} for further discussion), we expect that if $n$ is positive then $n \lessapprox (\log k)^2$. If $n$ is negative, then $n$ must be fairly close to $-k$ by similar reasoning.

We are able to prove Conjecture \ref{conj:b factor present} when $k$ is large and $\ell$ is small.

\begin{theorem}\label{thm:all terms, by^ell}
Let $k_0$ be a sufficiently large positive constant, and let $k \geq k_0$ be an integer. Let $\ell$ be an integer with $2 \leq \ell \leq (\log \log k)^{1/5}$. If there are nonzero integers $n,d,b,y$ with $d \geq 1$, $\textup{gcd}(n,d)=1$, and $P^+(b) \leq k$ such that \eqref{eq:eqn of study, all terms, by^ell} holds, then $P^+ \left( \prod_{i=0}^{k-1}(n+id)\right) \leq k$, $d \in \{1,2\}$, and $|n| \leq k(d+1)$. Furthermore, $(n,d)$ is one of the pairs in \eqref{eq:solutions with d=1,2 and b present}.
\end{theorem}

It is possible to generalize the situation yet further. It is not so difficult to find triples $(k,n,d)$ that ``almost'' yield a solution to \eqref{eq:eqn of study, all terms, by^ell}, but not quite. For instance, if $(k,n,d)=(6,-36,17)$, then
\begin{align*}
\prod_{i=0}^5 (n+id) = (-36)(-19)(-2)(15)(32)(49)=-2^8 \cdot 3^3 \cdot 5 \cdot 7^2 \cdot 19.
\end{align*}
The prime $19>k$ occurs to an odd power in the $i=1$ term, so the full product cannot be written as $by^2$ with $P^+(b) \leq k$. However, if we omit the $i=1$ term, then
\begin{align*}
\prod_{\substack{i=0 \\ i \neq 1}}^5 (n+id)  = 15\cdot (2^4 \cdot 3 \cdot 7)^2
\end{align*}
is of the form $by^2$ with $P^+(b) \leq k$. Furthermore, we even have that $y$ is divisible by a prime $>k$, which is not a phenomenon we observed above. These types of ``sporadic'' solutions with large values of $d$ are quite interesting, and become more common when one allows oneself to drop terms from the product. However, as $k$ increases, it appears that one eventually must have $d \in \{1,2\}$ again, even after dropping a term. For instance, $k=19$ appears to be the last time one gets a solution to
\begin{align*}
\prod_{\substack{i=0 \\ i \neq 1}}^{18} (n+id) = by^2
\end{align*}
with $d >2$ (take $d=3$ and $n \in \{-34,-26\}$).

Of course, once we allow ourselves to drop a single term from the product, there is no reason we cannot allow ourselves to drop more terms. If we take $(k,n,d)=(7,6,23)$, for example, and drop the $i=1$ and $i=2$ terms, then we obtain the interesting
\begin{align*}
\prod_{\substack{i=0 \\ i \neq 1,2}}^6 (n+id) = 6 \cdot 75 \cdot 98 \cdot 121 \cdot 144 = (2^3 \cdot 3^2 \cdot 5\cdot 7 \cdot 11)^2.
\end{align*}
Even so, the solutions with $d>2$ seem to disappear as $k$ increases; if we drop the $i=1$ and $i=2$ terms then the last such solutions seem to appear at $k=19$, and if we drop the $i=1$ and $i=3$ terms then the last such solutions seem to appear at $k=23$.

Accordingly, it is no surprise that a number of authors have investigated variants of \eqref{eq:eqn of study, all terms, by^ell} in which one is allowed to drop some number of terms (e.g. \cite{Bennett2023,DLS2018,Edis2019,Pratt2025,Saradha1997,SS2001,SN1996,ST1992}). For our contribution, we prove the following theorem.

\begin{theorem}\label{thm:main general thm}
Let $k_0$ be a sufficiently large positive constant, and let $k \geq k_0$ be an integer. Let $\ell$ be an integer with $2 \leq \ell \leq (\log \log k)^{1/5}$. Assume $S \subseteq \{0,1,2,\ldots,k-1\}$ with $|S| \geq k - 0.2 \frac{k}{\log k}$. If there are nonzero integers $n,d,b,y$ with $d \geq 1$, $\textup{gcd}(n,d)=1, P^+(b) \leq k$, and
\begin{align}\label{eq:main general eqn of study}
\prod_{i \in S} (n+id) = by^\ell,
\end{align}
then $P^+\left(\prod_{i \in S} (n+id) \right) \leq k$, $d \leq 3(\log k)^{3/2}$, and $|n| \leq k(d+1)$. If $|S| \geq k - \frac{k}{(\log k)^3}$, then $d \in \{1,2\}$.
\end{theorem}

\subsection{Some heuristic discussion}\label{subsec:heuristic}

It is not completely clear what analogue of Conjectures \ref{conj:main conj} and \ref{conj:b factor present} we should make for solutions to \eqref{eq:main general eqn of study}. Theorem \ref{thm:main general thm} allows one to remove $\approx \frac{k}{\log k}$ terms $n+id$ from the product over $0 \leq i \leq k-1$. In this regime, the only solutions to \eqref{eq:main general eqn of study} have $P^+(y) \leq k$. However, if one removes enough terms, then we expect \eqref{eq:main general eqn of study} to have solutions with $P^+(y) > k$.

To see this, consider the following situation, which is based on a heuristic analysis involving smooth numbers (see Granville's excellent survey \cite{Gra2008} for an introduction to smooth numbers). Let $k$ be large, let $\ell \geq 2$ be small, and let $p$ be a prime just larger than $k$. Set $n=p^\ell$. If $d \neq p$ is a prime with $d \lessapprox n/k$, then standard conjectures on smooth numbers in short intervals and arithmetic progressions (e.g. \cite[Section 4]{Gra2008}) suggest that
\begin{align*}
\sum_{\substack{n <m \leq n+(k-1)d \\ P^+(m) \leq k \\ m \equiv n \, (\text{mod }d)}} 1 \approx \rho(\ell) k,
\end{align*}
where $\rho$ is the Dickman--de Bruijn function (see \cite[p. 268]{Gra2008}). If we let $W \subseteq \{1,\ldots,k-1\}$ denote the set of those $i$ with $n+id = m \in (n,n+(k-1)d]$ and $P^+(m) \leq k$, then we ought to have $|W| \gtrapprox \rho(\ell) k$. Furthermore,
\begin{align*}
\prod_{i \in W \cup \{0\}}(n+id) = by^\ell
\end{align*}
for some positive integers $b$ and $y$ with $P^+(b) \leq k$ and $P^+(y) > k$, since $p \mid y$.

As a concrete example of these heuristics, if $k=40,n=41^2,d=29$, and
\begin{align*}
W = \{1,3,11,13,15,16,19,20,27,37,38,39\},
\end{align*}
then $\prod_{i \in W\cup \{0\}} (n+id) = by^2$ with
\begin{align*}
b &= 5 \cdot 7 \cdot 17 \cdot 19 \cdot 23 \cdot 31 \cdot 37, \\
y &= 2^{11} \cdot 3^5 \cdot 5^2 \cdot 7^2 \cdot 11^2 \cdot 13 \cdot 17 \cdot 19 \cdot 23 \cdot 41.
\end{align*}
Note that $|W| = 12$ and $\rho(\ell)\cdot k = \rho(2) \cdot 40 = 12.274\ldots$.

We leave it as an open problem to formulate a reasonable conjecture on solutions to \eqref{eq:main general eqn of study} that takes into account the sizes of $k,\ell,d,|S|$, and $P^+(y)$.

\begin{remark}
The results in this paper have been formalized in Lean (assuming some standard mathematical results). See our Github repository \cite{Pratt2026Lean} for further information.
\end{remark}

\subsection{Notation and conventions}

Throughout the paper, we generally assume that $k$ is larger than some sufficiently large (effectively computable) absolute constant. We always assume $\ell$ is an integer with $2 \leq \ell \leq (\log \log k)^{1/5}$. Note that $\ell$ need not be prime.

We use the usual asymptotic notations $o(\cdot),O(\cdot),\ll,\gg$. 

Given a nonzero integer $n$, we write $P^+(n)$ for the largest prime dividing $n$ if $|n| \neq 1$, and if $n = \pm 1$ then we define $P^+(n)=1$.

If $m$ and $q$ are coprime integers, we write $\overline{m} \pmod{q}$ for the multiplicative inverse of $m$ modulo $q$.

For a prime $p$ and a nonzero integer $n$, we write $v_p(n)$ for the $p$-adic valuation of $n$.

\subsection{Outline of the paper}

We give an outline of the arguments in Section \ref{sec:outline of arguments}. In Section \ref{sec:preliminaries and proofs of main thms}, we give some preliminary results that are important for all that follows in the paper (the reader is particularly encouraged to pay Lemmas \ref{lem:factorization of terms} and \ref{lem:throw away to get size control} special heed). We then state Proposition \ref{prop:main prop}, which is the main result from which we deduce our main theorems. We prove Theorem \ref{thm:main general thm} assuming Proposition \ref{prop:main prop}. We then prove Theorem \ref{thm:all terms, by^ell} using Theorem \ref{thm:main general thm}, and prove Theorem \ref{thm:explicit solutions} using Theorem \ref{thm:all terms, by^ell}.

In Section \ref{sec:important lemmas}, we prove two important lemmas, Lemma \ref{lem:remonds theorem} and Lemma \ref{lem:many m+1 sets div by large f}, and then state Proposition \ref{prop:pos dens from mass inc}. We then prove Proposition \ref{prop:main prop} using the two lemmas and Proposition \ref{prop:pos dens from mass inc}.

We prove Proposition \ref{prop:pos dens from mass inc} in Section \ref{sec:proof of last key prop}. Proposition \ref{prop:pos dens from mass inc} depends, in turn, on Proposition \ref{prop:few reps of alpha has large size}. Among other results, Proposition \ref{prop:few reps of alpha has large size} relies on Lemma \ref{lem:few alpha's have large multiplicity}.

\section{Outline of the arguments}\label{sec:outline of arguments}

The arguments in this paper have their roots in earlier work of the author \cite{Pratt2024} on rational points on Erd\H{o}s--Selfridge curves. In particular, an adaptation of the ``mass-increment argument'' introduced in \cite{Pratt2024} plays a vital role. However, in this paper, we do not assume familiarity with \cite{Pratt2024}. We give in this section an overview of the arguments in this paper, suppressing technical details for ease of exposition.

Theorems \ref{thm:explicit solutions} and \ref{thm:all terms, by^ell} are deduced from the more general Theorem \ref{thm:main general thm}, along with some additional number-theoretic inputs like results on special cases of \eqref{eq:eqn of study, no b} (see Lemma \ref{lem:Saradha square lemma}). 

The main task is therefore to prove Theorem \ref{thm:main general thm}, which concerns solutions to \eqref{eq:main general eqn of study}. It is not so hard to check that any solution to \eqref{eq:main general eqn of study} induces a factorization $n+id=a_ix_i^\ell$ for $i \in S$, where $a_i$ is only divisible by primes $\leq k$, and the integers $x_i$ are pairwise coprime and only divisible by primes $>k$ (Lemma \ref{lem:factorization of terms}). The main thrust of the method is to study the interplay between the $a_i$ and the $x_i$. In particular, Proposition \ref{prop:main prop} states that there are many $i \in S$ such that $|a_i| \leq k$ and $|x_i| = 1$. By a simple sort of bootstrapping argument, this implies that $|x_i|=1$ for every $i \in S$, which, in conjunction with results on primes in arithmetic progressions (Lemma \ref{lem:prime in AP}), implies Theorem \ref{thm:main general thm}.

The key ingredient in proving Proposition \ref{prop:main prop} is Proposition \ref{prop:pos dens from mass inc}, which says that the set of $a_i$ with $|a_i| \leq k$ has ``large'' cardinality. We prove Proposition \ref{prop:main prop} by contradiction. If there are not many $i \in S$ such that $|a_i| \leq k$ and $|x_i| =1$, then by Proposition \ref{prop:pos dens from mass inc} we can find many $|a_i| \leq k$ such that $|x_i| \neq 1$. Using the pigeonhole principle and other combinatorial ingredients (see Lemma \ref{lem:many m+1 sets div by large f}), we find there is some superelliptic curve of ``low height'' that has many rational points. (We know the rational points are distinct from each other because they involve the $x_i$, and $x_i$ with $|x_i| \neq 1$ are pairwise coprime and therefore distinct.) The superelliptic curve actually has ``too many'' rational points, contradicting quantitative results on the number of rational points on curves of genus at least two (Lemma \ref{lem:remonds theorem}).

Proposition \ref{prop:pos dens from mass inc} is the heart of the proof, and it is proved by an adaptation of the mass-increment argument introduced in \cite{Pratt2024}. We sketch the main ideas here briefly (see also \cite[Section 3]{Pratt2024} for a longer sketch).

The integers $a_i$ may be very large (we cannot initially rule out some of them are divisible by an arbitrarily large power of two, for instance), but after throwing away some of the $a_i$, we have a ``large'' set $I \subseteq S$ with
\begin{align}\label{eq:sketch up bound factorial}
\prod_{i \in I} |a_i| \leq k!
\end{align}
(see Lemma \ref{lem:throw away to get size control}). If $|a_i| \geq k$, then $|a_i|$ is distinct (Lemma \ref{lem:large ai are unique}), and an argument with Stirling's formula and \eqref{eq:sketch up bound factorial} implies there are some $i \in I$ with $|a_i| \leq k$. This raises the possibility of an iterative argument involving \eqref{eq:sketch up bound factorial}, provided we can show that
\begin{align*}
\prod_{\substack{i \in I \\ |a_i| \leq k}} |a_i|
\end{align*}
has some nontrivial size (the actual proof does not take an iterative approach, but it is helpful motivation to think of the proof in this way). We accomplish this by showing that the $a_i$ ``mostly'' have to be distinct from one another. The main result we need is Lemma \ref{lem:few alpha's have large multiplicity}, which says there are few $\alpha$ such that $\alpha = |a_i|$ for many different $i \in I$. The proof of Lemma \ref{lem:few alpha's have large multiplicity} relies on a variant of the pigeonholing argument described above, and the argument ultimately relies on the same bound for the number of rational points on a curve of genus at least two.

We close this section by mentioning some limitations of our approach. Recall the heuristic discussion in section \ref{subsec:heuristic} above, with $n=p^\ell$, and $S$ chosen so that the $n+id$ with $i \in W$ are all $k$-smooth. In terms of the factorizations $n+id=a_ix_i^\ell$, it is clear that $a_0=1,x_0=p$, and $a_i=n+id,x_i=1$ for all $i \in S$. Therefore, in this case, almost all of the $a_i$ are large, and the combinatorial ideas underlying the mass-increment argument fail. Therefore, new ideas will be needed to study \eqref{eq:main general eqn of study} in this situation.

\section{Preliminaries, and proofs of the main theorems}\label{sec:preliminaries and proofs of main thms}

We need a result on the structure of the terms $n+id$ in \eqref{eq:main general eqn of study}.

\begin{lemma}\label{lem:factorization of terms}
Let $k$ and $\ell$ be integers $\geq 2$. Let $S$ be a nonempty subset of $\{0,1,\ldots,k-1\}$. Let $n,d,b,y$ be integers with $bny \neq 0, d \geq 1$, $\textup{gcd}(n,d)=1$, and $P^+(b) \leq k$ such that \eqref{eq:main general eqn of study} holds. Then for $i \in S$ there are nonzero integers $a_i,x_i$ such that the following hold:
\begin{itemize}
\item $n+id$ may be factored uniquely as $n+id = a_ix_i^\ell$.
\item If $p$ is a prime with $p \mid a_i$, then $p \leq k$.
\item If $p$ is a prime with $p \mid x_i$, then $p > k$.
\item $\textup{gcd}(a_i,a_j) \mid (j-i)$ for all distinct $i,j \in S$.
\item $\textup{gcd}(a_i,d) = \gcd(x_i,d)=1$.
\item $\textup{gcd}(x_i,x_j)=1$ for all distinct $i,j \in S$.
\end{itemize}
\end{lemma}
\begin{proof}
We have $n+id \neq 0$ since $by \neq 0$. Let $i,j \in S$ be distinct. Then $\text{gcd}(n+id,n+jd) \mid (j-i)d$, and since $n$ and $d$ are coprime we have $\text{gcd}(n+id,n+jd) \mid (j-i)$. Note that $|j-i| \leq k-1$. In particular, this implies that if a prime $p$ divides two different terms $n+id$ and $n+jd$, then $p < k$.

We factor $n+id = a_i z_i$, where $a_i$ is only divisible by primes $\leq k$ and $z_i$ is only divisible by primes $> k$. We have $\text{gcd}(z_i,z_j)=1$ if $i \neq j$, so each $z_i$ is an $\ell$-th power by \eqref{eq:main general eqn of study} and the fact that $P^+(b) \leq k$. We accordingly write $z_i = x_i^\ell$. It now follows that $\text{gcd}(x_i,x_j) = 1$, and $\text{gcd}(a_i,a_j) \mid(j-i)$. For uniqueness, we always require $x_i$ to be positive, and $a_i$ has the sign of $n+id$.

Lastly, we note that since $n+id = a_ix_i^\ell$ and $\text{gcd}(n,d)=1$, we have $\text{gcd}(a_i,d)=1$ and $\gcd(x_i,d)=1$.
\end{proof}

Initially, we have little control over the sizes of the integers $a_i$ in Lemma \ref{lem:factorization of terms}. However, by throwing away a small number of elements of $S$, we can gain more control.

\begin{lemma}\label{lem:throw away to get size control}
Let $k\geq 10$ be an integer. Let $S \subseteq \{0,1,\ldots,k-1\}$ satisfy $|S| > \pi(k)$. For each $i \in S$, let $b_i$ be a nonzero integer divisible only by primes $\leq k$. Assume $\textup{gcd}(b_i,b_j) \mid (j-i)$ for distinct $i,j \in S$. Then, there exists $I \subseteq S$ with $|I| \geq |S|-\pi(k)$ such that
\begin{align*}
\prod_{i \in I} |b_i| \leq k!.
\end{align*}
\end{lemma}
\begin{proof}
For a prime $p \leq k$, let $i_p \in S$ be such that $v_p(b_{i_p})$ is maximal. We claim that $I = S \backslash \{i_p : p \leq k\}$ has the desired property. We trivially have $|I| \geq |S|-\pi(k)$ by the union bound.

It suffices to show $v_p(\prod_{i \in I} b_i) \leq v_p(k!)$. Fix a prime $p \leq k$. Then
\begin{align*}
v_p \left(\prod_{i \in I} b_i \right) &\leq v_p \left(\prod_{\substack{i \in S \\ i \neq i_p}} b_i \right) \leq v_p \left(\prod_{\substack{i \in S \\ i \neq i_p}} \text{gcd}(b_i,b_{i_p}) \right) \leq v_p \left(\prod_{\substack{i=0 \\ i \neq i_p}}^{k-1} |i-i_p| \right) \\
&= v_p(i_p! (k-1-i_p)!) \leq v_p ((k-1)!). \qedhere
\end{align*}
\end{proof}

By Lemmas \ref{lem:factorization of terms} and \ref{lem:throw away to get size control}, there exists $I \subseteq S$ with $|I| \geq |S|-\pi(k)$ such that
\begin{align}\label{eq:prod bd for ai with i in I}
\prod_{i \in I} |a_i| \leq k!.
\end{align}
Whenever we write $I$ in this paper, we shall mean this particular set $I$.

We require the following version of the prime number theorem in arithmetic progressions.

\begin{lemma}\label{lem:prime in AP}
Let $x \geq y \geq 5$ be real numbers. Let $d$ be a positive integer with $2\leq d \leq (\log x)^{1.99}$, and let $a$ be an integer coprime to $d$. Then, for any fixed $A>0$, we have
\begin{align*}
\sum_{\substack{y < p \leq x \\ p \equiv a (d)}} \log p = \frac{x-y}{\phi(d)} + O_A \left(\frac{x}{(\log x)^A} \right),
\end{align*}
where the implied constant is effectively computable.
\end{lemma}
\begin{proof}
If $d=2$, this is a consequence of the prime number theorem with classical error term, so we may assume $d \geq 3$. By the discussion in \cite[Chapter 20]{Dav2000}, we have
\begin{align}\label{eq:prime in AP class num}
\sum_{\substack{y < p \leq x \\ p \equiv a (d)}} \log p = \frac{x-y}{\phi(d)} + O \left(\frac{x}{\phi(d)} \exp \left(-c_1 \frac{\log x}{d^{1/2}(\log d)^2} \right) + x \exp(-c_2 (\log x)^{1/2}) \right),
\end{align}
where $c_1$ and $c_2$ are absolute, effectively computable positive constants. Since $d \leq (\log x)^{1.99}$, we obtain
\begin{align*}
\frac{\log x}{d^{1/2}(\log d)^2} \gg \frac{\log x}{(\log x)^{0.995}(\log \log x)^2}\gg (\log x)^{0.001}.
\end{align*}
Hence, for any fixed $A>0$, we have
\begin{align*}
\exp \left(-c_1 \frac{\log x}{d^{1/2}(\log d)^2} \right) + \exp(-c_2 (\log x)^{1/2}) \ll_A (\log x)^{-A},
\end{align*}
where the implied constant is effectively computable.
\end{proof}

\begin{remark}
The error term $O\left(\frac{x}{\phi(d)} \exp \left(-c_1 \frac{\log x}{d^{1/2}(\log d)^2} \right)\right)$ in the proof of Lemma \ref{lem:prime in AP} arises from using Dirichlet's effective class number bound on a potential Siegel zero \cite[(12) on p. 96]{Dav2000}. Since $d$ is small compared to $x$, this class number bound suffices to bound the error, and we do not need to appeal to Siegel's bound with ineffective constants \cite[Chapter 21]{Dav2000}.
\end{remark}

For the duration of the paper, unless explicitly stated otherwise, we fix a solution of \eqref{eq:main general eqn of study} satisyfing the hypotheses of Theorem \ref{thm:main general thm}. The factors $a_i,x_i$ are those of Lemma \ref{lem:factorization of terms}, and $I$ is a set obtained by Lemma \ref{lem:throw away to get size control}.

The following proposition is the main result upon which our theorems depend.

\begin{proposition}\label{prop:main prop}
Let $k$ be large and $2 \leq \ell \leq (\log \log k)^{1/5}$. Let $S \subseteq \{0,1,\ldots,k-1\}$ with $|S| \geq k-0.2\frac{k}{\log k}$. With notation as above, we have
\begin{align*}
\#\{i \in I : |a_i| \leq k \textup{ and } |x_i| = 1\} > \frac{k}{(\log k)^{3/2}}.
\end{align*}
\end{proposition}

\begin{proof}[Proof of Theorem \ref{thm:main general thm} assuming Proposition \ref{prop:main prop}]
By Proposition \ref{prop:main prop}, we can find $i_0,j_0 \in I$ such that $|a_{i_0}|\leq k, |a_{j_0}| \leq k, |x_{i_0}| = |x_{j_0}|=1$, and $|i_0-j_0| \geq \frac{k}{(\log k)^{3/2}} - 1$. Then
\begin{align*}
|i_0-j_0|d &= |(n+i_0d)-(n+j_0d)| \leq |n+i_0d| + |n+j_0d| = |a_{i_0}x_{i_0}^\ell| + |a_{j_0}x_{j_0}^\ell| \leq 2k,
\end{align*}
and therefore
\begin{align*}
d \leq \frac{2k}{\frac{k}{(\log k)^{3/2}}-1} \leq 3 (\log k)^{3/2},
\end{align*}
the last inequality holding for $k\geq 12$.

Note that
\begin{align*}
|n| &= |n+i_0d - i_0d| \leq |n+i_0d| + i_0d = |a_{i_0}x_{i_0}^\ell| + i_0 d \leq k + kd = k(d+1).
\end{align*}
Hence, for any $i \in S$, we have
\begin{align*}
|n+id| &\leq |n| + id \leq k(d+1) + kd \leq 6k(\log k)^{3/2} + k.
\end{align*}
If $|x_i| \neq 1$, then $|x_i| > k$, and then
\begin{align*}
k^2 &\leq |x_i|^\ell \leq |a_ix_i^\ell| = |n+id| \leq 6k(\log k)^{3/2} + k,
\end{align*}
and this is a contradiction for $k\geq 46$. It follows that $|x_i| = 1$ for every $i \in S$. Lemma \ref{lem:factorization of terms} then implies $P^+\left(\prod_{i \in S} (n+id) \right) \leq k$.

It remains to show that if $|S| \geq k - \frac{k}{(\log k)^3}$, then $d \in \{1,2\}$. We show we cannot have $d\geq 3$. 

Assume for contradiction that $3 \leq d \leq 3(\log k)^{3/2}$. We consider two cases, depending on the size of $n$. 

First, assume $n> - \frac{3k}{2}$. Then $n+id > -\frac{3k}{2} + id$, and $-\frac{3k}{2} + id>k$ if $i > \frac{5k}{2d}$. If $\frac{5k}{2d} <i \leq k-1$, then $n+\frac{5k}{2} < n+id \leq n+(k-1)d$. By Lemma \ref{lem:prime in AP}, we have
\begin{align*}
\sum_{\substack{n+\frac{5k}{2} < p \leq n+(k-1)d \\ p \equiv n (d)}}  \log p = \frac{kd(1 - \frac{5}{2d})}{\phi(d)} + O \left(\frac{k}{(\log k)^2}\right).
\end{align*}
The main term has size $\gg k$ since $d \geq 3$. Hence, there are $\gg \frac{k}{\log k}$ different $i$ with $\frac{5k}{2d} <i \leq k-1$ such that $n+id$ is a prime $>k$. There are at most $\frac{k}{(\log k)^3}$ different $i \not \in S$, so we conclude there is some $i \in S$ such that $n+id$ is a prime $>k$. This is a contradiction.

Now assume $n \leq -\frac{3k}{2}$. Write $m=|n|$, so that $m\geq \frac{3k}{2}$ and $n+id = -(m-id)$. By Lemma \ref{lem:prime in AP}, we have
\begin{align*}
\sum_{\substack{k < p \leq m \\ p \equiv m(d)}} \log p = \frac{m-k}{\phi(d)} + O \left(\frac{k}{(\log k)^2} \right).
\end{align*}
The main term is $\gg \frac{k}{d} \gg \frac{k}{(\log k)^{3/2}}$, so there are $\gg \frac{k}{(\log k)^{5/2}}$ different $i\geq 0$ such that $p=m-id$ is a prime $>k$. Note that $k<m-id$ if and only if $i < \frac{m-k}{d} \leq \frac{k(d+1)-k}{d}=k$. By arguing as before, we see there is some $i \in S$ such that $m-id$ is a prime $>k$, and this is a contradiction.
\end{proof}

We next deduce Theorem \ref{thm:all terms, by^ell} from Theorem \ref{thm:main general thm}. 

\begin{proof}[Proof of Theorem \ref{thm:all terms, by^ell}]
Assume we have a solution to \eqref{eq:eqn of study, all terms, by^ell}. By Theorem \ref{thm:main general thm}, we have $d \in \{1,2\}$ and $|n| \leq k(d+1) \leq 3k$. We also have that no term $n+id$ is divisible by a prime $> k$. To finish the proof, we must show that $(n,d)$ is one of the solutions in \eqref{eq:solutions with d=1,2 and b present}. Recall that we write $p$ for the least prime $>k$.

First, we assume that $n$ is positive. We claim that we must have $d=1$. Suppose for contradiction that $d=2$. It then follows from Lemma \ref{lem:prime in AP} (with $d=2$) and the fact that $n \ll k$ that one of the terms $n,n+2,\ldots,n+2(k-1)$ is a prime $>k$, and this is a contradiction. 

It follows that $d=1$. Now we claim that $n+k-1 < p$. If $n \geq p$, then every term $n+i$ exceeds $k$, and by the prime number theorem with classical error term, some term $n+i$ is a prime, a contradiction. Hence $n < p$. If $n+k-1 \geq p$, then some term $n+i=p$, and this is again a contradiction. Hence $n+k-1 < p$. All terms $n+i$ are positive and smaller than $p$, so every prime divisor is at most $k$. Therefore, the solutions $(n,d)$ with $n$ positive are exactly
\begin{align*}
(n,1), \ \ \ \ 1 \leq n \leq p-k.
\end{align*}

Next, we assume that $n$ is negative. We have two subcases to consider, depending on whether $d=1$ or $d=2$.

We first assume $d=1$. Since $by \neq 0$, we must have $n < -(k-1)$, so $n+(k-1) < 0$. Let $m = n+(k-1)$. Then
\begin{equation}\label{eq:switching trick}
\begin{split}
by^\ell &= \prod_{j=0}^{k-1}(n+j) = \prod_{j=0}^{k-1}(m-(k-1)+j) = (-1)^k \prod_{j=0}^{k-1} (-m+(k-1-j)) \\
&= (-1)^k\prod_{i=0}^{k-1} (-m+i).
\end{split}
\end{equation}
Since $-m > 0$, we may use the previous case to find that the solutions with $n$ negative and $d=1$ are
\begin{align*}
(n,1), \ \ \ \ -(p-1) \leq n \leq -k.
\end{align*}

Now assume that $n$ is negative and we have a solution with $d=2$. Since $\text{gcd}(n,d)=1$, we must have that $n$ is odd. We claim that $n > -2(k-1)$. Otherwise, if $n < -2(k-1)$, we can use the ``switching trick'' as in \eqref{eq:switching trick} to obtain a solution with $m$ positive and $d=2$, and we already saw this is not possible. We therefore have $-2(k-1) < n < 0$. This implies $n > -p$. If not, then $n \leq -p$, but then there is some $i$ such that $n+2i = -p$, and this contradicts the fact that no term is divisible by a prime $>k$. A bit of thought then reveals that all the solutions have
\begin{align*}
(n,2), \ \ \ \ -(p-2) \leq n \leq p-2-2(k-1), \ \ \ \ n \text{ odd}.
\end{align*}

Therefore, in any case, $(n,d)$ is one of the solutions in \eqref{eq:solutions with d=1,2 and b present}.
\end{proof}

We need the following result before embarking on the proof of Theorem \ref{thm:explicit solutions}.

\begin{lemma}\label{lem:Saradha square lemma}
Let $k\geq 2$ be an integer. There are no solutions to the equation
\begin{align*}
\prod_{i=0}^{k-1}(m+2i)=y^2
\end{align*}
in odd positive integers $m$ and $y$.
\end{lemma}
\begin{proof}
If $k\geq 3$, then the result follows from work of Saradha \cite[Corollary 1]{Saradha1997}. It therefore suffices to consider the case $k=2$.

Suppose for contradiction that we have $m(m+2) = y^2$ with $m,y$ positive and odd. Since $m$ and $m+2$ are coprime, we have $m=a^2,m+2=b^2$ for some positive integers $a < b$. Subtracting gives $2 = b^2-a^2 = (b+a)(b-a)$, so we must have $b-a = 1$ and $b+a=2$. Solving this pair of equations gives $a = \frac{1}{2}, b = \frac{3}{2}$, which is a contradiction.
\end{proof}

\begin{proof}[Proof of Theorem \ref{thm:explicit solutions}]
Assume we have a solution to \eqref{eq:eqn of study, no b} with $k$ sufficiently large and $2 \leq \ell \leq (\log \log k)^{1/5}$. By Theorem \ref{thm:all terms, by^ell} with $b=1$, we have $d \in \{1,2\}$ and $|n| \leq k(d+1) \leq 3k$. Furthermore, no term $n+id$ is divisible by a prime $> k$.

We show first that $\ell=2$. Let $p \in (\frac{k}{2},k)$ be a prime. Since $p \nmid d$, we have $p \mid(n+id)$ if and only if $i \equiv -n \overline{d} (p)$. Given the size of $p$, we see there is at least one $i$ such that $p \mid (n+id)$, and there are at most two such $i$. Then $p^\ell$ divides the left-hand side of \eqref{eq:eqn of study, no b}, and since $\text{gcd}(n+id,n+jd) < k$, we see that if $\ell \geq 3$ then there is some $i$ such that $p^2 \mid (n+id)$. But
\begin{align*}
\frac{k^2}{4} < p^2 \leq |n+id| \leq |n|+id \leq |n|+2i \leq 3k+2k = 5k,
\end{align*}
and this is a contradiction for $k\geq 21$. We conclude that $\ell=2$.

We next show we cannot have $d=1$. If $d=1$, then we have $\prod_{i=0}^{k-1} (n+i) = y^2$ for some $|n| \leq 2k$. The work of Erd\H{o}s and Selfridge \cite[Theorem 1]{ES1975} implies $n$ cannot be $\geq 0$. If $n$ is negative, then $n\leq -k$, and the change of variables $t=-n-(k-1)$ (recall \eqref{eq:switching trick}) yields
\begin{align*}
\prod_{j=0}^{k-1}(t+j) = (-1)^ky^2,
\end{align*}
where $1 \leq t \leq k+1$. There are no solutions if $k$ is odd (since the left-hand side is positive but the right-hand side is negative), and if $k$ is even there are again no solutions, by the work of Erd\H{o}s and Selfridge. Hence $d \neq 1$. (Actually, in the special case $\ell=2$, the full Erd\H{o}s--Selfridge result is not required, and one can appeal to earlier work by Erd\H{o}s \cite{Erd1939b} and Rigge \cite{Rigge1939}.)

We have shown that $d=2$. Therefore, we have an odd integer $n$ such that
\begin{align*}
\prod_{i=0}^{k-1} (n+2i) = y^2.
\end{align*}
If $n$ is positive, then we obtain a contradiction by applying Lemma \ref{lem:Saradha square lemma}, so we must have that $n$ is negative. If $n < -2(k-1)$, then we change variables $t=-n-2(k-1)>0$ to get
\begin{align*}
(-1)^k y^2 = \prod_{i=0}^{k-1}(t+2(k-1-i)) = \prod_{j=0}^{k-1} (t+2j).
\end{align*}
We have an immediate contradiction if $k$ is odd. If $k$ is even, then we get a contradiction from Lemma \ref{lem:Saradha square lemma}. Hence $n > -2(k-1)$.

We have shown that $-2(k-1) < n < 0$. Thus, some of the terms $n+2i$ are negative, and some are positive. Let $N,P$ denote the number of negative and positive terms, respectively. Clearly $N+P=k$. The negative terms are $n,n+2,\ldots,n+2(N-1)$, and the positive terms are $n+2N,\ldots,n+2(k-1)$. Since the integers $n+2i, 0 \leq i \leq k-1$, are consecutive odd integers, we see that $n+2(N-1) = -1$, and $n=-(2N-1)$. It follows that
\begin{align*}
y^2 = (-1)^N \prod_{i=0}^{N-1}(2i+1) \cdot \prod_{j=0}^{P-1}(2j+1).
\end{align*}
Clearly we must have $N$ is even, so
\begin{align}\label{eq:NP square}
y^2 = \prod_{i=0}^{N-1}(2i+1) \cdot \prod_{j=0}^{P-1}(2j+1).
\end{align}

We claim that \eqref{eq:NP square} has no solutions if $P \geq N+2$ or $N \geq P+2$. Indeed, suppose for contradiction that we have a solution to \eqref{eq:NP square} with $P = N+L, L \geq 2$. Then we have
\begin{align*}
y^2 = \prod_{i=0}^{N-1}(2i+1)^2 \cdot \prod_{j=N}^{N+L-1}(2j+1),
\end{align*}
so $\prod_{j=N}^{N+L-1}(2j+1)$ is equal to a square, but this is impossible by Lemma \ref{lem:Saradha square lemma}. We argue similarly in the case $N \geq P+2$.

We have proved that, if \eqref{eq:NP square} has a solution, then $P-1 \leq N \leq P+1$. We consider these three different cases in turn.

If $N = P+1$ then there are no solutions. Otherwise, we would have
\begin{align*}
y^2 = (2(N-1)+1)\prod_{j=0}^{P-1} (2j+1)^2,
\end{align*}
so $2(N-1)+1 = 2N-1$ is a square. But $N$ is even, so this is impossible modulo 4.

Assume $N=P$. Then $k=N+P=2N$ is divisible by four, and $n=-(k-1)$.

Lastly, assume $N=P-1$. It is easy to check that $2P-1$ must be a square. Then $k=N+P=(P-1)+P=2P-1$ is an odd square, and $n=-(2N-1) = -(2P-3) = -(k-2)$.
\end{proof}

\section{Important lemmas, and the proof of Proposition \ref{prop:main prop}}\label{sec:important lemmas}

We first introduce some lemmas we need for the proof of Proposition \ref{prop:main prop}.

The following result is a special case of work of R\'emond \cite{Rem2010} on the number of rational points on a curve. It is a quantitative version of Faltings' theorem \cite{Faltings1983} that a curve of genus at least two has at most finitely many rational points.

\begin{lemma}\label{lem:remonds theorem}
Let $\ell,m$ be integers with $\ell \geq 2, m \geq \max(2,7-\ell)$. Let $B,C_1,C_2,\ldots,C_m$ be nonzero integers, with $C_1,C_2,\ldots,C_m$ distinct. Let $F(X,Y)$ be the polynomial in two variables given by
\begin{align*}
F(X,Y) = BY^\ell - \prod_{j=1}^m (X+C_j),
\end{align*}
and let $H$ be an upper bound for the absolute value of the coefficients of $F$. Let $D = \max(\ell,m)$ denote the degree of $F$. Then the number of points $(X,Y) \in \mathbb{Q}^2$ with $F(X,Y) = 0$ is
\begin{align*}
\leq\exp(5^{D^4}(\log 3H)(\log \log 3H)).
\end{align*}
\end{lemma}
\begin{proof}
By the hypotheses on $\ell,m$, and the integers $C_j$, we see $F$ defines a smooth affine curve. The genus of the smooth projective completion is
\begin{align*}
\frac{(\ell-1)(m-1)-\gcd(\ell,m)+1}{2},
\end{align*}
which is $\geq 2$ if $m \geq \max(2,7-\ell)$ (see, e.g., \cite[Exercise A.4.6]{HS2000} for the genus calculation). Hence, by Faltings' theorem, $F$ has only finitely many rational zeros. The result then follows from \cite[Th\'eor\`eme 1.1]{Rem2010}.
\end{proof}

\begin{remark}\label{rmk:upper bound on ell not optimal given remond}
We work with the assumption $\ell \leq (\log \log k)^{1/5}$ in this paper, but the true limit imposed by Lemma \ref{lem:remonds theorem} is $\ell \leq \delta (\log \log k)^{1/4}$ for some suitably small constant $\delta > 0$.
\end{remark}

\begin{remark}\label{rmk:better bds on num of rat pts}
Lemma \ref{lem:remonds theorem} is the reason our main theorems are restricted to $\ell \leq (\log \log k)^{1/5}$. This strict condition arises from the doubly-exponential dependence of the bound in Lemma \ref{lem:remonds theorem} on the degree $D$. R\'emond's theorem \cite[Th\'eor\`eme 1.1]{Rem2010} works for general curves, but one should be able to do much better in our situation by exploiting the superelliptic structure of $F$. 

When $\ell$ is prime, one ought to be able to use work of Poonen--Schaefer \cite{PS1997} to show that the rank $r$ of the Jacobian of the curve defined by $F$ satisfies $r \ll (\ell m)^{O(1)}\log (3H)$. When combined with recent work on quantitative bounds for Mordell's conjecture due to Yu, Yuan, and Zhou \cite{YYZ2026}, one should obtain a bound of $\ll \exp((\ell m)^{O(1)}\log (3H))$ for the number of rational points on the superelliptic curve. This would then allow us to prove versions of our main theorems with the condition $\ell \leq (\log \log k)^{1/5}$ replaced by $\ell \leq (\log k)^\delta$, where $\delta>0$ is a small constant.
\end{remark}

The second result is a generalization of \cite[Proposition 4.3]{Pratt2024}, which goes back to Erd\H{o}s \cite[Lemma 1]{Erdos1939}. The result says that if we have a set of integers $\leq k$ of sufficiently large cardinality, then there are many tuples of integers in the set all divisible by a large common factor.

\begin{lemma}\label{lem:many m+1 sets div by large f}
Let $k,m$ be positive integers. Let $1 \leq b_1< b_2 < \cdots < b_{r-1} < b_r \leq k$ be positive integers. Let $c \in (0,1]$ be such that $r \geq ck$, and assume $\frac{m}{c} \leq \frac{\log k}{\log \log k}$. Then there are $\geq \frac{ck}{6(m+1)}$ pairwise disjoint sets $\{b_{i_0},b_{i_1},\ldots,b_{i_{m}}\}$ with $m+1$ elements each such that, for every such set, there is an integer $f > \frac{c}{3}\exp(-2.77\frac{m}{c})k$ with $f \mid b_{i_j}$ for $0 \leq j \leq m$.
\end{lemma}
\begin{proof}
Fix a parameter $\eta > 0$ to be chosen later. Let $\mathcal{F}$ denote the set of positive integers $n$ in $(\eta k , k]$  such that every proper divisor $f$ of $n$ satisfies $f \leq \eta k$. Every integer in $(\eta k,k]$ is divisible by some $f \in \mathcal{F}$.

Let $\mathcal{B}$ denote the set of $b_i$ with $\eta k < b_i \leq k$. Note that $|\mathcal{B}| \geq r - \eta k$. Define a function $g: \mathcal{B} \rightarrow \mathcal{F}$ given by $g(b) = f$, where $f$ is the largest element of $\mathcal{F}$ that divides $b$. By double-counting, we have
\begin{align*}
|\mathcal{B}| &= \sum_{b \in \mathcal{B}} \sum_{\substack{f \in \mathcal{F} \\ f = g(b)}} 1 = \sum_{f \in \mathcal{F}} \sum_{\substack{b \in \mathcal{B} \\ g(b) = f}} 1 = \sum_{f \in \mathcal{F}} \#U(f),
\end{align*}
say, so that $U(f) = \{b \in \mathcal{B} : g(b) = f\}$. Note that the sets $U(f)$ are disjoint for different $f \in \mathcal{F}$. We split $f \in \mathcal{F}$ depending on $\#U(f)$, and trivial estimation then yields
\begin{align*}
|\mathcal{B}| -m|\mathcal{F}| &\leq \sum_{\substack{f \in \mathcal{F} \\ \#U(f) \geq m+1}} \#U(f).
\end{align*}

If $\#U(f) \geq m+1$, then there are $\left\lfloor \frac{\#U(f)}{m+1} \right\rfloor \geq 1$ disjoint subsets $\{b_{i_0},b_{i_1},\ldots,b_{i_{m}}\}$ of $U(f)$, each with $m+1$ elements, such that $f \mid b_{i_j}$ for $0 \leq j \leq m$. The total number of such $(m+1)$-element sets is then
\begin{align*}
&\geq \sum_{\substack{f \in \mathcal{F} \\ \#U(f) \geq m+1}}\left\lfloor \frac{\#U(f)}{m+1} \right\rfloor \geq \frac{1}{2(m+1)} \sum_{\substack{f \in \mathcal{F} \\ \#U(f) \geq m+1}} \#U(f) \geq \frac{|\mathcal{B}|-m|\mathcal{F}|}{2(m+1)}.
\end{align*}

In order to finish the proof, we need an upper bound on $|\mathcal{F}|$. Let $A>0$ be a parameter to choose shortly, and set $\delta = A\eta$. Then
\begin{align*}
|\mathcal{F}| &\leq \delta k + \#\{f \in \mathcal{F} : \delta k < f \leq k\}.
\end{align*}
Let $f \in \mathcal{F}$ satisfy $f > \delta k$, and let $p$ be the least prime divisor of $f$. Since $f \in \mathcal{F}$, we have $\frac{f}{p} \leq \eta k$, and therefore $p > A$. Hence, if $f \in \mathcal{F}$ and $f > \delta k$, then $P^-(f) > A$. It follows that
\begin{align*}
|\mathcal{F}| &\leq A \eta k + \sum_{\substack{n \leq k \\ P^-(n) > A}} 1.
\end{align*}
By the fundamental lemma of sieve theory (see, e.g., \cite[Lemma 6.8]{Opera}), if $A \leq k^{1/100}$, say, then we have
\begin{align*}
\sum_{\substack{n \leq k \\ P^-(n) > A}} 1 \leq k(1+O(e^{-s})) \prod_{p \leq A} \left(1 - \frac{1}{p} \right),
\end{align*}
where $s = \frac{1}{10} \frac{\log k}{\log A}$. It follows that
\begin{align*}
|\mathcal{B}| - m|\mathcal{F}| &\geq k \left(c - \eta (mA+1) - m(1+O(e^{-s})) \prod_{p \leq A} \left(1 - \frac{1}{p}\right) \right).
\end{align*}
If $A \geq 3$, then by work of Rosser and Schoenfeld \cite[Theorem 7]{RS1962} we have
\begin{align*}
\prod_{p \leq A} \left(1 - \frac{1}{p}\right) \leq \frac{0.794}{\log A}.
\end{align*}
We set $A = \exp(2.4 \frac{m}{c})$. Note that $A \geq \exp(2.4) \geq 11$, and $A \leq k^{O(1/\log \log k)} \leq k^{1/100}$ by our upper bound on $\frac{m}{c}$. Since $s \rightarrow \infty$ as $k \rightarrow \infty$, we have
\begin{align*}
m(1+O(e^{-s})) \prod_{p \leq A} \left(1 - \frac{1}{p}\right) \leq \frac{0.8 m}{\log A} = \frac{c}{3}.
\end{align*}
Observe that
\begin{align*}
mA+1 = m\exp\left(2.4 \frac{m}{c}\right)+1 \leq \frac{m}{c}\exp\left(2.4 \frac{m}{c}\right)+1 \leq \exp\left(2.77 \frac{m}{c}\right).
\end{align*}
If we set $\eta = \frac{c}{3} \exp(-2.77 \frac{m}{c})$, then
\begin{align*}
\frac{|\mathcal{B}| - m|\mathcal{F}|}{2(m+1)} \geq \frac{ck}{6(m+1)},
\end{align*}
as desired.
\end{proof}

We introduce an important parameter $m_\ell$, which is a positive integer that depends on $\ell$. We define
\begin{align}\label{eq:defn of m_ell}
m_\ell = 
\begin{cases}
6, &\ell \in \{2,3\}, \\
\ell, &\ell \geq 4.
\end{cases}
\end{align}
The significance of $m_\ell$ will become clear in the proofs of Proposition \ref{prop:main prop} and Lemma \ref{lem:few alpha's have large multiplicity}. Specifically, $m_\ell$ is chosen so that certain curves have genus at least two, and so that $\ell \mid m_\ell$.

\begin{proposition}\label{prop:pos dens from mass inc}
Assume $\#\{i \in I : |a_i| \leq k \textup{ and } |x_i| = 1\} \leq \frac{k}{(\log k)^{3/2}}$. If $2 \leq \ell \leq (\log \log k)^{1/5}$, then
\begin{align*}
\#\{1 \leq \alpha \leq k : \textup{there exists }i \in I \textup{ with } |a_i| = \alpha\} \geq \frac{k}{1000 m_\ell \log m_\ell}.
\end{align*}
\end{proposition}

\begin{proof}[Proof of Proposition \ref{prop:main prop} assuming Proposition \ref{prop:pos dens from mass inc}]
Assume by way of contradiction that 
\begin{align*}
\#\{i \in I : |a_i| \leq k \text{ and } |x_i| = 1\} \leq \frac{k}{(\log k)^{3/2}}.
\end{align*}
Write $m = m_\ell$, and define
\begin{align*}
\mathcal{A} := \{1 \leq \alpha \leq k : \text{there exists }i \in I \text{ with } |a_i| = \alpha \text{ and } |x_i| \neq 1\}.
\end{align*}
Then, by Proposition \ref{prop:pos dens from mass inc}, we have
\begin{align*}
|\mathcal{A}| \geq \frac{k}{2000 m \log m} = c_m \cdot k.
\end{align*}

We apply Lemma \ref{lem:many m+1 sets div by large f} to the set $\mathcal{A}$ with $c = c_m$. (Note that $\frac{m}{c_m} = 2000m^2\log m = o(\frac{\log k}{\log \log k})$.) If we define 
\begin{align*}
\eta_m = \frac{c_m}{3}\cdot \exp \left(-2.77 \frac{m}{c_m} \right),
\end{align*}
then we deduce there are
\begin{align*}
\geq \frac{c_m k}{6(m+1)}
\end{align*}
disjoint sets $\{a_{i_0},a_{i_1},\ldots,a_{i_m}\}$ with $m+1$ elements each, with a common divisor $f > \eta_m k$.

Consider one such set $\{a_{i_0},a_{i_1},\ldots,a_{i_m}\}$, where we may assume $i_0 < i_1 < \cdots < i_m$. For $1 \leq j \leq m$, we have
\begin{align*}
(i_j-i_0)d = a_{i_j}x_{i_j}^\ell - a_{i_0} x_{i_0}^\ell.
\end{align*}
We have $f \mid a_{i_j}$ and $f\mid a_{i_0}$, so $f \mid (i_j-i_0) d$. However, we have $\text{gcd}(f,d)=1$, so $f \mid (i_j-i_0)$. Hence, dividing by $f$ and rearranging, we obtain
\begin{align*}
\frac{a_{i_j}}{f}x_{i_j}^\ell = \frac{a_{i_0}}{f}x_{i_0}^\ell + \frac{i_j-i_0}{f} \cdot d.
\end{align*}
We multiply these expressions together to obtain
\begin{align*}
(x_{i_1}\cdots x_{i_{m_\ell}})^\ell \cdot \prod_{j=1}^{m} \frac{a_{i_j}}{f} = \prod_{j=1}^{m} \left(\frac{a_{i_0}}{f}x_{i_0}^\ell + \frac{i_j-i_0}{f} \cdot d \right).
\end{align*}
We divide both sides by $d^{m_\ell}$, and use the fact that $\ell \mid m$. If we set $X = \frac{a_{i_0}}{f}\frac{x_{i_0}^\ell}{d}, Y = \frac{x_{i_1}\cdots x_{i_m}}{d^{m/\ell}}$, then $(X,Y)$ is a rational point on a curve of the form
\begin{align}\label{eq: BY^ell = prod X curves}
B Y^\ell = \prod_{j=1}^m (X+C_j),
\end{align}
where $B$ is a nonzero integer with $|B| \leq \eta_m^{-m}$, and $1 \leq C_1 < C_2 < \cdots < C_m$ are distinct positive integers with $C_j \leq \eta_m^{-1}$. The total number of curves of this form is $\leq 3 \eta_m^{-2m}$. A curve of this form has the absolute value of its coefficients bounded by $H \leq (2\eta_m^{-1})^m$.

Note that, for different sets $\{a_{i_0},a_{i_1},\ldots,a_{i_{m_\ell}}\}$, the associated rational points $(X,Y)$ are distinct. This is because all the $x_i$ under consideration have $|x_i| \neq 1$, and the $x_i$ are pairwise coprime by Lemma \ref{lem:factorization of terms}.

By the pigeonhole principle, some curve of the form \eqref{eq: BY^ell = prod X curves} has
\begin{align}\label{eq:low bd num rat pts}
\gg k \exp(-O\left(m^3 \log m \right) ) = \exp((1-o(1))\log k)
\end{align}
rational points.

We apply Lemma \ref{lem:remonds theorem} with $D = m, H = (2\eta_m^{-1})^m$ to see that any curve of the form \eqref{eq: BY^ell = prod X curves} has 
\begin{align}\label{eq:up bd num rat pts}
\ll \exp\left(5^{m^4} \cdot O(m^3 (\log m)^2)\right) = \exp(o(\log k))
\end{align}
rational points. Clearly, equations \eqref{eq:low bd num rat pts} and \eqref{eq:up bd num rat pts} contradict each other.
\end{proof}

\section{Proof of Proposition \ref{prop:pos dens from mass inc}}\label{sec:proof of last key prop}

This section is devoted to proving Proposition \ref{prop:pos dens from mass inc}. Before stating the results on which Proposition \ref{prop:pos dens from mass inc} depends, we need to introduce some notation and definitions.

Given a positive integer $\alpha \leq k$, define the set $T(\alpha) = \{i \in I : |a_i| = \alpha\}$. Note that the sets $T(\alpha)$ are disjoint for distinct $\alpha$. We have
\begin{align}\label{eq:num of i in I related to sum of T(alpha)}
\#\{i \in I : |a_i| \leq k\} = \sum_{\substack{i \in I \\ |a_i| \leq k}} 1 = \sum_{i \in I} \sum_{\substack{\alpha \leq k \\ \alpha = |a_i|}} 1 = \sum_{\alpha \leq k} \sum_{\substack{i \in I \\ |a_i| = \alpha}}1 = \sum_{\alpha \leq k} \#T(\alpha).
\end{align}
Given a positive integer $m$, we split the sum over $\alpha$ according to the size of $T(\alpha)$:
\begin{align}\label{eq:split sum over T into leq m and gtr m}
\sum_{\alpha \leq k} \#T(\alpha) = \sum_{\substack{\alpha \leq k \\ \#T(\alpha) \leq m}} \#T(\alpha) + \sum_{\substack{\alpha \leq k \\ \#T(\alpha) >m}} \#T(\alpha).
\end{align}

Recall the definition of $m_\ell$ in \eqref{eq:defn of m_ell}. The following result immediately yields Proposition \ref{prop:pos dens from mass inc}.

\begin{proposition}\label{prop:few reps of alpha has large size}
Let $m=m_\ell$, and define $\rho_m = \frac{1}{1000\log m}$. Assume $\#\{i \in I : |a_i| \leq k \text{ and } |x_i| = 1\} \leq \frac{k}{(\log k)^{3/2}}$. Then
\begin{align*}
\sum_{\substack{\alpha \leq k \\ \#T(\alpha) \leq m}} \#T(\alpha)> \rho_m \cdot k.
\end{align*}
\end{proposition}
\begin{proof}[Proof of Proposition \ref{prop:pos dens from mass inc} assuming Proposition \ref{prop:few reps of alpha has large size}]
Observe that
\begin{align*}
m\sum_{\substack{\alpha \leq k \\ 1\leq \#T(\alpha) \leq m}}1 \geq \sum_{\substack{\alpha \leq k \\ 1 \leq \#T(\alpha) \leq m}} \#T(\alpha) > \rho_m \cdot k.
\end{align*}
We finish the proof by noting
\begin{align*}
\#\{1 \leq \alpha \leq k : \text{there exists }i \in I \text{ with } |a_i| = \alpha\} &\geq \sum_{\substack{\alpha \leq k \\ 1\leq \#T(\alpha) \leq m}}1. \qedhere
\end{align*}
\end{proof}

\subsection{Proof of Proposition \ref{prop:few reps of alpha has large size}}

The first lemma shows that if $i \in S$ with $|a_i| \geq k$, then no other $j \in S$ with $j \neq i$ can have $|a_j| = |a_i|$.

\begin{lemma}\label{lem:large ai are unique}
Assume $i \in S$ with $|a_i| \geq k$. Then for all $j \in S$ with $j \neq i$, we have $|a_j| \neq |a_i|$.
\end{lemma}
\begin{proof}
We proceed by contraposition. Hence, assume we have distinct $i,j \in S$ with $|a_i| = |a_j|$. We may write $a_i = \alpha,a_j = \epsilon \alpha$, where $\epsilon \in \{\pm 1\}$. Since $n+id=a_ix_i^\ell$ and $n+jd = a_jx_j^\ell$, we see that
\begin{align*}
(i-j)d = (n+id)-(n+jd) = a_ix_i^\ell - a_jx_j^\ell = \alpha(x_i^\ell - \epsilon x_j^\ell),
\end{align*}
so $\alpha \mid (i-j)d$. We have $\text{gcd}(\alpha,d)=1$ by Lemma \ref{lem:factorization of terms}, so $\alpha \mid (i-j)$, and therefore $|a_i| = |\alpha| \leq |i-j| \leq k-1$.
\end{proof}

We can now obtain a lower bound on the left-hand side of \eqref{eq:num of i in I related to sum of T(alpha)}.

\begin{lemma}\label{lem:there are some i in I with ai < k}
Assume $|S| \geq k- 0.2 \frac{k}{\log k}$. If $k$ is sufficiently large, then 
\begin{align*}
\#\{i \in I : |a_i| \leq k\} \geq 0.1 \frac{k}{\log k}.
\end{align*}
\end{lemma}
\begin{proof}
Let $R$ denote the number of $i \in I$ with $|a_i| > k$. By Lemma \ref{lem:large ai are unique}, each such $|a_i|$ is unique. Hence, by \eqref{eq:prod bd for ai with i in I}, we have
\begin{align*}
\prod_{j=0}^{R-1}(k+j) \leq \prod_{\substack{i \in I \\ |a_i| > k}} |a_i| \leq k!.
\end{align*}
We take logarithms of both sides and then estimate to get
\begin{align*}
(R+k)\log(R+k)-k\log k-R \leq k\log k - k + O(\log k).
\end{align*}
Since $R < k$, we write $R = (1-\frac{\beta}{\log k}) k$, where $\beta > 0$. We plug in, rearrange, and cancel common terms to get
\begin{align*}
\left(2- \frac{\beta}{\log k} \right) \log \left(2- \frac{\beta}{\log k} \right)\cdot  k \leq \beta k + O(\log k).
\end{align*}
If $\beta \leq 1.35$, then we get a contradiction, since the left-hand side would then be $\geq (1-o(1)) (2 \log 2) k \geq (1.385\ldots )k$. We deduce that $\beta > 1.35$. Since
\begin{align*}
|S|-\pi(k) \leq |I| = \#\{i \in I : |a_i| < k\} + \#\{i \in I : |a_i| \geq k\} = \#\{i \in I : |a_i| < k\} + R,
\end{align*}
it follows that
\begin{align*}
\#\{i \in I : |a_i| \leq k\} &\geq |S| - \pi(k) - R \geq k - 0.2 \frac{k}{\log k} - \left(1-\frac{1.35}{\log k}\right) k \\
&= 1.15\frac{k}{\log k} - \pi(k)\geq 0.1 \frac{k}{\log k},
\end{align*}
where we have used the prime number theorem in the last inequality.
\end{proof}

We need some understanding of those $\alpha < k$ with $\#T(\alpha) > m_\ell$.

\begin{lemma}\label{lem:few alpha's have large multiplicity}
Assume $\#\{i \in I : |a_i| \leq k \text{ and } |x_i| = 1\} \leq \frac{k}{(\log k)^{3/2}}$. Then
\begin{align*}
\sum_{\substack{\alpha < k \\ \#T(\alpha) > m_\ell}} \#T(\alpha) \leq 0.049 \frac{k}{\log k}.
\end{align*}
\end{lemma}

\begin{lemma}\label{lem:many distinct gives lower bound increment}
Let $m$ be a positive integer with $4 \leq m \leq \sqrt{k}$. Let $\rho_m = \frac{1}{1000 \log m}$. Assume
\begin{align*}
\sum_{\substack{\alpha \leq k \\ \#T(\alpha) \leq m}} \#T(\alpha) \geq \delta_0 \frac{k}{\log k}
\end{align*}
for some $\delta_0 \in [0.001,\rho_m \cdot \log k]$. Then
\begin{align*}
\#\{i \in I : |a_i| \leq k\} &\geq (\delta_0 + 0.1) \frac{k}{\log k}.
\end{align*}
\end{lemma}
\begin{proof}
Let us define $\mathcal{A}_m = \{\alpha \leq k : 1\leq \#T(\alpha) \leq m\}$. If we define $\mathcal{B} = \{i \in I : |a_i| \in \mathcal{A}_m\}$, then by assumption we have $|\mathcal{B}|=J \geq \delta_0 \frac{k}{\log k}$. Let us temporarily relabel the $a_i$'s so that $|a_1| \leq |a_2| \leq \cdots \leq |a_J|$ are the $|a_i|$ with $i \in \mathcal{B}$. Note that
\begin{align*}
\prod_{\substack{i \in I \\ |a_i| \leq k}} |a_i| \geq \prod_{j=1}^J |a_j| \geq 1^m \cdot \prod_{j=m+1}^J |a_j| \geq 1^m \cdot 2^m \cdot \prod_{j=2m+1}^J |a_j| \geq \cdots.
\end{align*}
In general, we have
\begin{align*}
\prod_{\substack{i \in I \\ |a_i| \leq k}} |a_i| &\geq \prod_{t=1}^{\lfloor J/m\rfloor} t^m = \left(\left\lfloor \frac{J}{m}\right\rfloor !\right)^m \geq \left(\left\lfloor \frac{\delta_0 k}{m\log k}\right\rfloor !\right)^m.
\end{align*}
Let $R$ denote the number of $i \in I$ such that $|a_i| > k$. By \eqref{eq:prod bd for ai with i in I} and Lemma \ref{lem:large ai are unique}, we then have
\begin{align*}
k! &\geq \prod_{i \in I} |a_i| = \prod_{\substack{i \in I \\ |a_i| \leq k}} |a_i| \cdot \prod_{\substack{i \in I \\ |a_i| > k}} |a_i| \geq \left(\left\lfloor \frac{\delta_0 k}{m\log k}\right\rfloor !\right)^m\cdot \prod_{j=0}^{R-1} (k+j).
\end{align*}
Let $y = \delta_0/\log k$, so $0 < y \leq \rho_m$. Taking logarithms and using Stirling's formula, we find
\begin{align*}
k \log k - k + O(m\log k) &\geq \delta_0 k -yk(\log(m/y)+1)+ (R+k)\log(R+k)-k\log k-R.
\end{align*}

The function $y(\log(m/y)+1)$ is increasing for $0 < y \leq m$, and since $y \leq \rho_m \leq 1 < m$, we have
\begin{align*}
y(\log(m/y)+1) \leq \rho_m (\log(m/\rho_m)+1) = \frac{1}{1000}\left(1 + \frac{\log \log m}{\log m} + \frac{\log(1000)+1}{\log m} \right).
\end{align*}
This latter function is a decreasing function of $m$, so we have
\begin{align*}
y(\log(m/y)+1) \leq \frac{1}{1000} \left(1 + \frac{\log \log 4}{\log 4} + \frac{\log(1000)+1}{\log 4} \right) \leq 0.007.
\end{align*}

It follows that
\begin{align*}
2k \log k - k + O(m \log k) \geq k(\delta_0 - 0.007) + (R+k)\log(R+k)-R.
\end{align*}
This implies $R < k$ since $m$ is small, so we may write $R = (1-\frac{\beta}{\log k})k$ with $\beta > 0$. We insert this expression for $R$, cancel some terms, and divide through by $k$ to obtain
\begin{align*}
\beta + O \left(\frac{m\log k}{k} \right) \geq \delta_0 - 0.007 + \left(2 - \frac{\beta}{\log k}\right) \log \left(2 - \frac{\beta}{\log k} \right).
\end{align*}
If we assume for contradiction that $\beta \leq \delta_0 + 1.37$, then we find the right-hand side is
\begin{align*}
\geq \delta_0 - 0.007 + (2-\rho_m-o(1))\log(2-\rho_m-o(1))\geq \delta_0 + 1.378,
\end{align*}
since $\rho_m \leq \frac{1}{1000\log 4}$. This yields a contradiction.

Hence, $\beta \geq \delta_0 + 1.37$. We find that
\begin{align*}
\#\{i \in I : |a_i| \leq k\} &\geq |S|-\pi(k)-R \geq k - 0.2 \frac{k}{\log k} - \pi(k) - \left(1 - \frac{\delta_0 + 1.37}{\log k}\right)k \\
&= (\delta_0 + 1.17)\frac{k}{\log k} - \pi(k) \geq (\delta_0 + 0.1) \frac{k}{\log k},
\end{align*}
where we have used the prime number theorem in the last step.
\end{proof}

\begin{proof}[Proof of Proposition \ref{prop:few reps of alpha has large size} assuming Lemma \ref{lem:few alpha's have large multiplicity}]
Write $m=m_\ell$. By Lemma \ref{lem:there are some i in I with ai < k} and \eqref{eq:num of i in I related to sum of T(alpha)}, we have
\begin{align*}
N = \sum_{\alpha \leq k} \#T(\alpha) \geq 0.1 \frac{k}{\log k}.
\end{align*}
We split the sum depending on the size of $T(\alpha)$, so that $N=J+H$, where
\begin{align*}
J &=\sum_{\substack{\alpha \leq k \\ \#T(\alpha) \leq m}} \#T(\alpha), \\
H &= \sum_{\substack{\alpha \leq k \\ \#T(\alpha) > m}} \#T(\alpha)
\end{align*}
Lemma \ref{lem:few alpha's have large multiplicity} gives $H \leq 0.049 \frac{k}{\log k}$, so $J \geq 0.051 \frac{k}{\log k}$.

Assume for contradiction that $J \leq \rho_m \cdot k$. If we set $\delta_0 = J \frac{\log k}{k}$, then $0.051 \leq \delta_0 \leq \rho_m \cdot \log k$ and $J \geq \delta_0 \frac{k}{\log k}$. It then follows from Lemma \ref{lem:many distinct gives lower bound increment} that
\begin{align*}
N \geq J + 0.1 \frac{k}{\log k}.
\end{align*}
However, 
\begin{align*}
N = J+H \leq J + 0.049 \frac{k}{\log k}
\end{align*}
by Lemma \ref{lem:few alpha's have large multiplicity}, which is a contradiction. We deduce that $J > \rho_m \cdot k$.
\end{proof}

\subsection{Proof of Lemma \ref{lem:few alpha's have large multiplicity}}

We need an upper bound on $\#T(\alpha)$.
\begin{lemma}\label{lem:bd on size of T alpha set}
Let $\alpha \leq k$ be a positive integer. Then $\#T(\alpha) \leq \frac{k}{\alpha} + 1$.
\end{lemma}
\begin{proof}
Let $i_1< i_2 < \cdots < i_J$ be the elements of $T(\alpha)$. We may assume $J \geq 2$, otherwise the desired bound is trivial. For any $j$ with $1 \leq j \leq J-1$, there are signs $\epsilon_j,\epsilon_{j+1} \in \{\pm 1\}$ such that
\begin{align*}
n+i_jd = \epsilon_j \alpha x_{i_j}^\ell, \ \ \ \ n+i_{j+1}d = \epsilon_{j+1}\alpha x_{i_{j+1}}^\ell.
\end{align*}
We therefore have $(i_{j+1}-i_j) d = \alpha (\epsilon_{j+1}x_{i_{j+1}}^\ell - \epsilon_jx_{i_j}^\ell)$, so $\alpha \mid (i_{j+1}-i_j) d$. Lemma \ref{lem:factorization of terms} implies $\alpha \mid (i_{j+1}-i_j)$, so $i_{j+1} \geq i_j + \alpha$. We deduce by induction that $i_j \geq (j-1)\alpha$ for $j\geq 1$, and since $i_J \leq k$ we have $J \leq \frac{k}{\alpha} + 1$, as desired.
\end{proof}

\begin{proof}[Proof of Lemma \ref{lem:few alpha's have large multiplicity}]
We write $m = m_\ell$ throughout. Assume for contradiction that $\#\{i \in I : |a_i| \leq k \text{ and } |x_i| = 1\} \leq \frac{k}{(\log k)^{3/2}}$ and
\begin{align*}
\sum_{\substack{\alpha < k \\ \#T(\alpha) > m}} \#T(\alpha) > 0.049 \frac{k}{\log k}.
\end{align*}
Let us say $\alpha$ is \emph{good} if we have $|x_i| \neq 1$ for all $i \in T(\alpha)$, and $\alpha$ is \emph{bad} if there is some $i \in T(\alpha)$ with $|x_i| = 1$. Observe that the number of bad $\alpha$ is $\leq \frac{k}{(\log k)^{3/2}}$.

We claim there are at most two $i \in T(\alpha)$ with $|a_i|=\alpha$ and $|x_i|=1$. If there are at least three, then there exist $i_1 < i_2 < i_3$ with
\begin{align*}
n+i_1d = a_{i_1}x_{i_1}^\ell, \ \ \ \ n+i_2d = a_{i_2}x_{i_2}^\ell, \ \ \ \ n+i_3d = a_{i_3}x_{i_3}^\ell,
\end{align*}
$|a_{i_j}| = \alpha$ for $j \in \{1,2,3\}$, and $|x_{i_j}|=1$ for $j \in \{1,2,3\}$. Hence, there are signs $\epsilon_j \in \{\pm 1\}$ for $j \in \{1,2,3\}$ with
\begin{align*}
n+i_1d = \epsilon_1 \alpha, \ \ \ \ n+i_2d = \epsilon_2 \alpha, \ \ \ \ n+i_3d = \epsilon_3 \alpha.
\end{align*}
By the pigeonhole principle, two of $n+i_1d,n+i_2d,n+i_3d$ are equal, and this is a contradiction since $d\neq 0$.

Given $\alpha$, define a set $V(\alpha)\subseteq T(\alpha)$ by $V(\alpha) = \{i \in I : |a_i| = \alpha \text{ and } |x_i| \neq 1\}$. By the claim above, we have $\#T(\alpha) \leq \#V(\alpha) + 2$. It follows that
\begin{align*}
0.049 \frac{k}{\log k} &\leq \sum_{\substack{\alpha \text{ good}\\ \#V(\alpha) > m}} \#V(\alpha) + \sum_{\substack{\alpha \text{ bad}\\ \#T(\alpha) \leq m+2}} \#T(\alpha)+\sum_{\substack{\alpha \text{ bad}\\ \#V(\alpha) > m}} (\#V(\alpha)+2) \\
&\leq \sum_{\substack{\alpha \leq k \\ \#V(\alpha) > m}} \#V(\alpha) + O \left(k \frac{m}{(\log k)^{3/2}} \right).
\end{align*}
Since $m \ll \ell\leq (\log k)^{1/4}$, say, we deduce
\begin{align*}
\frac{k}{\log k} \ll\sum_{\substack{\alpha \leq k \\ \#V(\alpha) > m}} \#V(\alpha).
\end{align*}
By dyadic pigeonholing, there is some $N$ with $1 \ll N \ll k$ such that
\begin{align}\label{eq:low bound on dyadic sum alpha with T alpha}
\frac{k}{(\log k)^2} \ll \sum_{\substack{N/2 \leq \alpha < N \\ \#V(\alpha) > m}} \#V(\alpha).
\end{align}
At this point we split into cases depending on the size of $N$.

First, assume $N > \frac{k}{(\log k)^3}$. Since $\#V(\alpha) \leq \#T(\alpha) \ll \frac{k}{\alpha} \ll \frac{k}{N}$ by Lemma \ref{lem:bd on size of T alpha set}, we have
\begin{align*}
\frac{N}{(\log k)^2} \ll \sum_{\substack{N/2 \leq \alpha < N \\ \#V(\alpha) > m}} 1.
\end{align*}
Hence, there are many $\alpha$ of size $\asymp N$ with $\#V(\alpha) > m$. Consider such an $\alpha$. By assumption, there exist $i_0 < i_1 < \cdots < i_m$ in $V(\alpha)$ with $|a_{i_j}| = \alpha$ for $0 \leq j \leq m$. We write $n+i_j d = \epsilon_j \alpha x_{i_j}^\ell$ with $\epsilon_j \in \{\pm 1\}$, and note that
\begin{align*}
\epsilon_j x_{i_j}^\ell = \epsilon_0 x_{i_0}^\ell + \frac{i_j - i_0}{\alpha} d
\end{align*}
for $1 \leq j \leq m$. Note that $\frac{i_j - i_0}{\alpha}$ is a positive integer (recall the proof of Lemma \ref{lem:bd on size of T alpha set}). We multiply together these $m$ different expressions to obtain
\begin{align*}
\epsilon_0^{m} \epsilon_1 \cdots \epsilon_{m} (x_{i_1}\cdots x_{i_{m}})^{\ell} &= \prod_{j=1}^{m} \left(x_{i_0}^\ell + \epsilon_0\frac{i_j-i_0}{\alpha}d\right).
\end{align*}

A key point is that $m$, the number of terms on the right-hand side, is divisible by $\ell$. Therefore, if we divide both sides by $d^{m}$, we obtain
\begin{align*}
\epsilon_0^{m} \epsilon_1 \cdots \epsilon_{m} \left(\frac{x_{i_1}\cdots x_{i_{m}}}{d^{m/\ell}}\right)^{\ell} &= \prod_{j=1}^{m} \left(\frac{x_{i_0}^\ell}{d} + \epsilon_0\frac{i_j-i_0}{\alpha}\right).
\end{align*}
We set 
\begin{align}
X = \frac{x_{i_0}^\ell}{d}, \ \ \ \ Y = \frac{x_{i_1}\cdots x_{i_{m}}}{d^{m/\ell}}
\end{align}
to see that $(X,Y)$ is a rational point on a curve of the form
\begin{align}\label{eq:form of curves from pigeonholing}
\varepsilon_1 Y^\ell = \prod_{j=1}^{m_\ell} (X+\varepsilon_2 B_j),
\end{align}
where $\varepsilon_1,\varepsilon_2 \in \{\pm 1\}$, and $1 \leq B_1 < B_2 < \cdots < B_m$ are distinct positive integers with each $B_j \leq 2(\log k)^3$. The total number of possible curves of this form is $\leq 4 \cdot 2^m (\log k)^{3m}$. 

Another point to note is that the rational points $(X,Y)$ arising in this way are distinct for distinct choices of $\alpha$. This is because the sets $V(\alpha)$ are disjoint, $i \in V(\alpha)$ implies $|x_i| \neq 1$, and the integers $x_i$ are pairwise coprime to each other by Lemma \ref{lem:factorization of terms}. We have at least $\gg \frac{N}{(\log k)^2} \gg \frac{k}{(\log k)^5}$ different rational points, so by the pigeonhole principle there is some curve of the form \eqref{eq:form of curves from pigeonholing} with at least
\begin{align*}
\gg \frac{k}{2^m (\log k)^{3m+5}} = \exp((1-o(1))\log k)
\end{align*}
rational points. However, by Lemma \ref{lem:remonds theorem}, the curve can have at most
\begin{align*}
\leq \exp(5^{m^4} (\log 3H)(\log \log 3H))
\end{align*}
rational points, where $H$ is an upper bound for the absolute value of the coefficients of the polynomial on the right-hand side of \eqref{eq:form of curves from pigeonholing}. A simple calculation shows
\begin{align*}
H \leq (1+2(\log k)^3)^m,
\end{align*}
and we find the curve has at most $\exp(o(\log k))$ distinct rational points, which is a contradiction.

Now consider the case in which $N \leq \frac{k}{(\log k)^3}$. Given a small constant $\epsilon_0 > 0$, the contribution of those $\alpha \asymp N$ to the right-hand side of \eqref{eq:low bound on dyadic sum alpha with T alpha} with $\#V(\alpha) \leq \epsilon_0 \frac{k}{N(\log k)^2}$ is $O(\epsilon_0 \frac{k}{(\log k)^2})$. Hence, if we take $\epsilon_0$ to be a sufficiently small constant, we have
\begin{align*}
\frac{k}{(\log k)^2} &\ll \sum_{\substack{N/2 \leq \alpha < N \\ \#V(\alpha) \gg kN^{-1}(\log k)^{-2}}} \#V(\alpha).
\end{align*}
Note that $\frac{k}{N(\log k)^2} \geq \log k$ is much larger than $m \ll \ell$. By Lemma \ref{lem:bd on size of T alpha set} we have $\#V(\alpha) \ll k/N$, and this implies
\begin{align*}
\frac{N}{(\log k)^2} \ll \sum_{\substack{N/2 \leq \alpha < N \\ \#V(\alpha) \gg kN^{-1}(\log k)^{-2}}} 1.
\end{align*}
The idea now is to get many different rational points from each $\alpha$, rather than just one rational point per $\alpha$ as in the case $N>\frac{k}{(\log k)^3}$.

Consider some $\alpha$ with $N/2 \leq \alpha < N$ and $\#V(\alpha) \gg \frac{k}{N(\log k)^2}$. The set $V(\alpha)$ is a subset of the interval $[0,k-1]$, and we cover $V(\alpha)$ by breaking this long interval into shorter intervals $J_r$. In particular, let $\Delta = N(\log k)^3$, and for nonnegative integers $r$ define $J_r = [r\Delta,(r+1)\Delta)$. We therefore have
\begin{align*}
\frac{k}{N(\log k)^2} &\ll \#V(\alpha) \ll \sum_{0 \leq r \leq k \Delta^{-1}} |V(\alpha) \cap J_r| \\ 
&= \sum_{\substack{0 \leq r \leq k \Delta^{-1} \\ |V(\alpha) \cap J_r| \leq m}} |V(\alpha) \cap J_r| + \sum_{\substack{0 \leq r \leq k \Delta^{-1} \\ |V(\alpha) \cap J_r| > m}} |V(\alpha) \cap J_r|.
\end{align*}
By trivial estimation we then have
\begin{align*}
\frac{k}{N(\log k)^2} &\ll \sum_{\substack{0 \leq r \leq k \Delta^{-1} \\ |V(\alpha) \cap J_r| > m}} |V(\alpha) \cap J_r|.
\end{align*}

We claim that $|V(\alpha) \cap J_r| \ll (\log k)^3$. Let $|V(\alpha) \cap J_r| = L$, where we may assume $L \geq 2$. Divide $J_r$ into $L-1$ subintervals, each of length $\frac{\Delta}{L-1}$. By the pigeonhole principle, some subinterval contains at least two distinct elements $i,j \in V(\alpha)$. By the familiar arguments, we have $\alpha \mid (i-j)$, so
\begin{align*}
\frac{N}{2} \leq \alpha \leq |i-j| \leq \frac{\Delta}{L-1} = \frac{N(\log k)^3}{L-1},
\end{align*}
and this gives the claim.

Applying the claim, we have
\begin{align*}
\frac{k}{N(\log k)^5} &\ll \sum_{\substack{0 \leq r \leq k \Delta^{-1} \\ |V(\alpha) \cap J_r| > m}} 1.
\end{align*}
Let $r$ be such that $|V(\alpha) \cap J_r| > m$. We now argue essentially identically to the argument above in the case $N > \frac{k}{(\log k)^3}$. We therefore obtain a rational point on a curve of the form \eqref{eq:form of curves from pigeonholing}. There are $\gg \frac{N}{(\log k)^2}$ choices of $\alpha \in [N/2,N)$ with $\# V(\alpha) \gg \frac{k}{N(\log k)^2}$, and, for each such $\alpha$, we get $\gg \frac{k}{N(\log k)^5}$ rational points on curves of the form \eqref{eq:form of curves from pigeonholing}. Since all the rational points are distinct from one another, this gives a total of
\begin{align*}
\gg \frac{N}{(\log k)^2} \cdot \frac{k}{N(\log k)^5} \gg \frac{k}{(\log k)^7}
\end{align*}
rational points, each of which is on some curve of the form \eqref{eq:form of curves from pigeonholing}. By the pigeonhole principle, there is at least one curve of this form with $\gg k 2^{-m}(\log k)^{-3m-7}$ rational points on it. This contradicts Lemma \ref{lem:remonds theorem}.
\end{proof}

\section*{Acknowledgments}

The author is partially supported by the National Science Foundation (DMS-2418328) and the Simons Foundation (MPS-TSM-00007959). 

ChatGPT 5.6 Sol was used for some literature search, translation, LaTeX assistance, proofreading, and for ``refereeing'' earlier versions of this manuscript. ChatGPT detected a small error in an earlier version of the proof of Lemma \ref{lem:many distinct gives lower bound increment} and suggested a fix. ChatGPT also suggested a slight simplification to an earlier version of the proof of Proposition \ref{prop:few reps of alpha has large size}. All the substantial mathematical ideas in the paper are due to the author.

The author used ChatGPT 6.1 Sol to formalize the results of this paper in Lean.

\bibliographystyle{plain}
\bibliography{references}

\end{document}